\documentclass[12pt,leqno]{amsart}
\usepackage[latin9]{inputenc}
\usepackage{verbatim}
\usepackage{amsthm}
\usepackage{amstext}
\usepackage{amssymb}
\usepackage[dvipdfmx]{graphicx}
\usepackage{xcolor}
\usepackage{soul}

\usepackage{mathrsfs}
\usepackage{amscd}

\makeatletter
\numberwithin{equation}{section}

\usepackage{amsthm}\usepackage{amscd}\usepackage{bm}

\usepackage[dvipdfmx,pdfstartview=FitH,pdfborderstyle={/S/B/W 1}]{hyperref}

\theoremstyle{plain}
\newtheorem{main theorem}{Main Theorem}
\newtheorem{theorem}{Theorem}[section]
\newtheorem{lemma}[theorem]{Lemma}

\newtheorem{proposition}[theorem]{Proposition}
\newtheorem{claim}[theorem]{Claim}

\theoremstyle{definition}
\newtheorem{definition}[theorem]{Definition}
\newtheorem{remark}[theorem]{Remark}

\newcommand{\mdim}{\mathrm{mdim}}
\newcommand{\vol}{\mathrm{vol}}

\makeatother

\begin{document}

\title[Remark on the local nature of metric mean dimension]{Remark on the local nature of metric mean dimension}

\author[Masaki Tsukamoto]{Masaki Tsukamoto}



\thanks{M. Tsukamoto was supported by JSPS KAKENHI 18K03275.}

\maketitle

\begin{abstract}
Metric mean dimension is a metric invariant of dynamical systems.
It is a dynamical analogue of Minkowski dimension of metric spaces.
We explain that old ideas of Bowen (1972) can be used for clarifying the local 
nature of metric mean dimension.
We also explain the generalization to $\mathbb{R}^D$-actions and an illustrating example.
\end{abstract}

\section{Introduction} \label{section: introduction}

\subsection{Background}  \label{subsection: background}

The purpose of this paper is to explain that old ideas of Bowen \cite{Bowen}
can be used in the context of \textit{mean dimension theory} 
\cite{Gromov, Lindenstrauss--Weiss, Lindenstrauss}.

A pair $(X, T)$ is called a \textbf{dynamical system} if 
$X$ is a compact metrizable space and $T:X\to X$ is a homeomorphism.
The mean dimension of $(X, T)$ (denoted by $\mdim(X, T)$)
is a topological invariant measuring the number of parameters \textit{per iterate}
for describing the orbits of $(X, T)$.
It was first introduced by Gromov \cite{Gromov}.

One of the big difficulties in the study of mean dimension is that 
it is often very hard to prove an upper bound on mean dimension. 
So far, the most successful way for proving upper bounds on $\mdim(X, T)$ is 
to use \textit{metric mean dimension}.
Metric mean dimension is a dynamical analogue of Minkowski dimension introduced by 
Lindenstrauss--Weiss \cite{Lindenstrauss--Weiss}.
It is not a topological invariant but a metric-dependent quantity.

Let $d$ be a metric (distance function) on $X$ compatible with the topology.
For $\varepsilon>0$ we denote by $S(X, \varepsilon)$ \textit{the entropy of $(X,T)$ at the scale $\varepsilon$}.
(The definition will be given in the next subsection.)
We define the \textbf{upper and lower metric mean dimensions} of $(X, T, d)$ by 
\begin{equation}  \label{eq: definition of metric mean dimension}
    \overline{\mdim}_M(X, T, d) := \limsup_{\varepsilon\to 0} \frac{S(X, \varepsilon)}{\log(1/\varepsilon)}, \quad 
     \underline{\mdim}_M(X, T, d)  := \liminf_{\varepsilon\to 0} \frac{S(X, \varepsilon)}{\log(1/\varepsilon)}. 
\end{equation}
They bounds the mean dimension $\mdim(X, T)$ from above \cite[Theorem 4.2]{Lindenstrauss--Weiss}:
\begin{equation} \label{eq: metric mean dimension bounds mean dimension}
   \mdim(X, T) \leq \underline{\mdim}_M(X, T, d) \leq \overline{\mdim}_M(X, T, d). 
\end{equation}
This is a dynamical analogue of the fact that Minkowski dimension bounds topological 
dimension. It is conjectured that there always exists a metric $d$ for which 
the equalities hold in (\ref{eq: metric mean dimension bounds mean dimension}).

The inequality (\ref{eq: metric mean dimension bounds mean dimension})
is a very useful result because it is often much easier to prove an upper bound on 
$\underline{\mdim}_M(X, T, d)$ than to (directly) prove an upper bound on $\mdim(X, T)$.
This idea was first successfully used in \cite{Tsukamoto Yang--Mills, Tsukamoto Brody}
for proving tight upper bounds on mean dimension of certain dynamical systems coming from geometric analysis. 

Why is it (relatively) easy to prove an upper bound on metric mean dimension?
The main reason is that we can calculate metric mean dimension by using only some \textit{local}
information.
(We will explain this more precisely in the next subsection.)
On the contrary, mean dimension has a highly global nature.
So far, we do not know any direct method to calculate it by using local information.
  
The purpose of this paper is to explain that
we can understand the local nature of metric mean dimension
more clearly by using ideas of Bowen \cite{Bowen}.

\subsection{Main result}  \label{subsection: main result}

Let $(X, d)$ be a compact metric space.
Let $\varepsilon>0$ and $K$ a subset of $X$.
A subset $S$ of $X$ is said to be an \textbf{$\varepsilon$-spanning set} of $K$
(or \textbf{$(d, \varepsilon)$-spanning set} of $K$ when we need to clarify what metric is used) 
if for every $x\in K$ there exists $y\in S$ satisfying $d(x, y) \leq \varepsilon$.
We denote by $\#(K, d, \varepsilon)$ the minimum cardinarity of such $S$.

Let $(X, T)$ be a dynamical system with a metric $d$ compatible with the topology.
For a subset $\Omega$ of $\mathbb{Z}$, we define a metric $d_\Omega$ on $X$ by 
\[ d_{\Omega}(x, y) := \sup_{a\in \Omega} d(T^a x, T^a y). \]
If $\Omega$ is a finite subset, then $d_\Omega$ is also compatible with the given topology.
However, if $\Omega$ is infinite (in particular, if $\Omega = \mathbb{Z}$), the metric $d_\Omega$
often defines a topology different from the originally given one.
(It is often the case that $(X, d_\Omega)$ is noncompact if $\Omega$ is infinite.)
For a natural number $N$, we also denote 
\[ d_N(x, y) := d_{\{0,1,2,\dots, N-1\}}(x, y) = \max_{0\leq n <N} d(T^n x, T^n y). \]

For $\varepsilon>0$ and a subset $K$ of $X$ we define 
\[  S(K, \varepsilon) := \limsup_{N\to \infty} \frac{\log \#(K, d_N, \varepsilon)}{N}. \]
This quantity is called the \textbf{entropy of $K$ at the scale $\varepsilon$}.
We define the upper and lower metric mean dimensions of $(X, T, d)$ as in 
(\ref{eq: definition of metric mean dimension}):
\[  \overline{\mdim}_M(X, T, d) := \limsup_{\varepsilon\to 0} \frac{S(X, \varepsilon)}{\log(1/\varepsilon)}, \quad 
     \underline{\mdim}_M(X, T, d)  := \liminf_{\varepsilon\to 0} \frac{S(X, \varepsilon)}{\log(1/\varepsilon)}. \]

Now we have defined metric mean dimension.
It is not difficult to see that its definition has a local nature:
If $X= K_1\cup K_2 \cup \dots \cup K_n$ then for any $\varepsilon>0$ and natural number $N$
\[ \#(X, d_N, \varepsilon) \leq \sum_{i=1}^n \#(K_i, d_N, \varepsilon). \]
This means that we can easily decompose the problem of estimating $\#(X, d_N, \varepsilon)$ into local pieces.
Then 
\[  S(X, \varepsilon) = \max_{1\leq i\leq n} S(K_i, \varepsilon) \]
and hence 
\[ \overline{\mdim}_M(X, T, d) = \limsup_{\varepsilon\to 0} 
   \frac{\max_{1\leq i\leq n}S(K_i, \varepsilon)}{\log(1/\varepsilon)}, \quad 
     \underline{\mdim}_M(X, T, d)  = \liminf_{\varepsilon\to 0} 
  \frac{\max_{1\leq i\leq n}S(K_i, \varepsilon)}{\log(1/\varepsilon)}. \]
So metric mean dimension can be calculated by local information (i.e. the entropy of each piece $K_i$ at the scale $\varepsilon$).
This is easy to see and not surprising. 
The main result below shows that indeed we can calculate metric mean dimension 
by using much more local quantity.

Let $\delta>0$ and $\Omega \subset \mathbb{Z}$.
For $x\in X$, we define $B_\delta(x, d_{\Omega})$ as the closed $\delta$-ball around $x$ with respect to 
$d_\Omega$:
\[  B_\delta(x, d_\Omega) := \{y\in X|\, d_\Omega(x, y) \leq \delta\}. \]
The following is our main result.

\begin{theorem}  \label{theorem: main result}
For any positive number $\delta$
   \begin{equation*}
       \begin{split}
      \overline{\mdim}_M(X, T, d) &= \limsup_{\varepsilon\to 0} 
   \frac{\sup_{x\in X}S\left(B_\delta(x, d_\mathbb{Z}), \varepsilon\right)}{\log(1/\varepsilon)}, \\
     \underline{\mdim}_M(X, T, d) & = \liminf_{\varepsilon\to 0} 
   \frac{\sup_{x\in X}S\left(B_\delta(x,d_{\mathbb{Z}}), \varepsilon\right)}{\log(1/\varepsilon)}.
 \end{split} 
   \end{equation*}
\end{theorem}

This is surprising because the set 
\[ B_\delta(x, d_{\mathbb{Z}}) = \{y\in X|\, d_\mathbb{Z}(x, y) \leq \delta\} = 
    \{y\in X|\, d(T^n x, T^n y) \leq \delta  \text{ for all integers $n$}\}  \]
is a very small subset.
It is usually \textit{not} a neighborhood of the point $x$.
For example, if $(X, T)$ is the full-shift on $\{0,1\}^\mathbb{Z}$ with some metric $d$
(or, more generally, if $(X, T)$ is \textit{expansive}), then 
the set $B_\delta(x, d_\mathbb{Z})$ is just one-point $\{x\}$ for sufficiently small $\delta$.
Nevertheless, we can calculate metric mean dimension 
only by studying such small subsets.

We prove Theorem \ref{theorem: main result} in the next section.
The proof is just a slight modification of the argument given by Bowen \cite{Bowen}.
Essentially, Theorem \ref{theorem: main result} is a mere simple corollary of \cite{Bowen}.
However Theorem \ref{theorem: main result}
seems to be useful in the future study\footnote{What I have in my mind is the problem of estimating 
mean dimension of various geometric examples in \cite{Gromov}.
For such examples, probably we can study the sets $B_\delta(x, d_\mathbb{Z})$ by some 
deformation theory techniques. See also \S \ref{section: example}.} of mean dimension theory, 
so the author thinks that it is worth publishing.

\begin{remark}
Bowen \cite{Bowen} considered the quantity 
\[  h^*_{T, \mathrm{homeo}}(\delta) 
  := \sup_{x\in X} \left\{\lim_{\varepsilon\to 0} S\left(B_\delta(x, d_{\mathbb{Z}}), \varepsilon\right)\right\}. \]
A dynamical system $(X, T)$ is said to be \textbf{$h$-expansive} if $h^*_{T, \mathrm{homeo}}(\delta) =0$
for some $\delta>0$.
He studied several important consequences of this condition. 
In the study of entropy theory of \cite{Bowen}, the entropy of $B_\delta(x, d_{\mathbb{Z}})$ is more or less a
``remainder term''. On the other hand, Theorem \ref{theorem: main result} shows that 
$S\left(B_\delta(x, d_{\mathbb{Z}})\right)$ becomes a ``main term'' in the study of metric mean dimension.
Therefore, interestingly, the viewpoints here and in \cite{Bowen} are almost opposite
although the proof of Theorem \ref{theorem: main result} is very close to \cite{Bowen}.
\end{remark}

\section{Proof of Theorem \ref{theorem: main result}}  \label{section: proof of main result}

Throughout this section, we assume that $(X, T)$ is a dynamical system with a metric $d$.
In this section, for two integers $a\leq b$, we write $[a, b]$ as the set of \textit{integers} $x$ with $a\leq x\leq b$.
(This convention is used only in this section. In the later sections, $[a, b]$ means the set of \textit{real numbers} $x$
with $a\leq x\leq b$.)

The next lemma is \cite[Lemma 2.1]{Bowen}.

\begin{lemma}[\cite{Bowen}]    \label{lemma: coding each piece}
Let $\varepsilon >0$ and $F\subset X$. Let 
\[  0=t_0 < t_1< t_2< \dots < t_{r-1} < t_r = n \]
be a sequence of integers. Then 
\[  \#(F, d_n, \varepsilon) \leq \prod_{i=0}^{r-1} \#\left(T^{t_i} F, \, d_{t_{i+1}-t_i},\,  \frac{\varepsilon}{2}\right). \]
\end{lemma}

\begin{proof}
Let $E_i$ be a $\left(d_{t_{i+1}-t_i}, \frac{\varepsilon}{2}\right)$-spanning set of $T^{t_i}(F)$ of minimum cardinality.
For each $(x_0, x_1, \dots, x_{r-1})\in E_0\times E_1\times \dots\times E_{r-1}$, consider 
\[ V(x_0, x_1, \dots, x_{r-1}) :=
 \left\{y\in F\middle|\, d_{t_{i+1}-t_i}\left(T^{t_i}y, x_i\right) \leq \frac{\varepsilon}{2} \text{ for all $0\leq i\leq r-1$}\right\}. \]
These form a covering of $F$. The diameter of $V(x_0, \dots, x_{r-1})$ with respect to $d_n$ is less than or equal to $\varepsilon$.
We take one point from each non-empty $V(x_0, x_1, \dots, x_{r-1})$.  
They form a $(d_n, \varepsilon)$-spanning set of $F$.
\end{proof}

The next proposition is a key result.
This is a very small modification of \cite[Proposition 2.2]{Bowen}.

\begin{proposition}  \label{prop: Bowen ball}
Let $\varepsilon$ and $\delta$ be positive numbers.
Set
\[ a := \sup_{x\in X} S\left(B_\delta(x,d_{\mathbb{Z}}), \frac{\varepsilon}{4}\right). \]
For any positive number $\beta$ we have 
\[  \sup_{x\in X} \#\left(B_\delta(x, d_n), d_n, \varepsilon\right) \leq C \cdot e^{(a+\beta)n} \quad 
    \text{for all integers $n$}, \]
where $C = C(\varepsilon, \delta, \beta)$ is a positive number depending on $\varepsilon, \delta, \beta$.
\end{proposition}

\begin{proof}
The proof is almost identical to \cite[Proposition 2.2]{Bowen}.
But we provide a full proof for the completeness.
For each $y\in X$, since we know that 
\[ \lim_{N\to \infty} \frac{1}{N} \log \#\left(B_\delta(y,d_{\mathbb{Z}}), d_N,\frac{\varepsilon}{4}\right) \leq a < a+ \beta, \]
there exist a natural number $m(y)$ and a subset $E(y)\subset X$ which is a $\left(d_{m(y)}, \frac{\varepsilon}{4}\right)$-spanning set 
of $B_\delta(y,d_{\mathbb{Z}})$ with $|E(y)| < e^{(a+\beta) m(y)}$.
We define an open set $U(y)$ by 
\[ U(y) := \left\{x\in X\middle|\, \exists z\in E(y): \> d_{m(y)} (x, z) < \frac{\varepsilon}{2}\right\}. \]
We have $B_\delta(y,d_{\mathbb{Z}}) \subset U(y)$ and 
\begin{equation} \label{eq: span of U(y)}
   \#\left(U(y), d_{m(y)}, \frac{\varepsilon}{2}\right) \leq |E(y)| < e^{(a+\beta) m(y)}. 
\end{equation}
We have 
\[  B_\delta(y, d_\mathbb{Z}) = \bigcap_{N=1}^\infty B_\delta\left(y, d_{[-N,N]}\right), \]
and each $B_\delta\left(y, d_{[-N,N]}\right)$ is closed.
Since $X$ is compact, there exists $N(y)>0$ satisfying 
\[  B_\delta\left(y, d_{[-N(y), N(y)]}\right) \subset U(y). \]
If we take a sufficiently small $\eta>0$ (depending on $y$) then we also have
\[  B_{\delta+\eta}\left(y, d_{[-N(y),N(y)]}\right) \subset U(y).  \] 
Set $V(y) := \{x\in X|\, d_{[-N(y), N(y)]}(x, y) < \eta\}$. This is an open neighborhood of $y$ and 
\[  \forall v\in V(y):  \quad B_\delta\left(v, d_{[-N(y),N(y)]}\right)  \subset U(y). \]
We choose $y_1, \dots, y_s\in X$ with $X = V(y_1)\cup \dots \cup V(y_s)$.
Set 
\[ N:= \max\left(N(y_1), \dots, N(y_s), m(y_1), \dots, m(y_s)\right) + 1. \]

Let $n>2N$ and $x\in X$. 
For each $t\in [N, n-N]$ we can choose $y_i$ with $T^t x\in V(y_i)$ and then 
\[  T^t\left(B_\delta(x, d_n)\right) \subset B_\delta\left(T^t x, d_{[-N(y_i), N(y_i)]}\right) \subset U(y_i). \]
From (\ref{eq: span of U(y)})
\[ \#\left( T^t\left(B_\delta(x, d_n)\right), d_{m(y_i)}, \frac{\varepsilon}{2}\right) \leq 
     \#\left(U(y_i), d_{m(y_i)}, \frac{\varepsilon}{2}\right) < e^{(a+\beta) m(y_i)}. \]

We will inductively choose integers $0=t_0< t_1< \dots< t_{r-1} < t_r=n$ and points $y_{i_1}, y_{i_2}, \dots, y_{i_{r-1}}$ with 
$T^{t_k} x\in V(y_{i_k})$ $(1\leq k\leq r-1)$.
First we set $t_1 = N$ and choose $y_{i_1}$ with $T^N x\in V(y_{i_1})$.
Next, suppose we have defined $t_1, \dots, t_k (< n)$ and $y_{i_1}, \dots, y_{i_k}$.
If $t_k>n-N$ then we set $r=k+1$ and $t_r = n$. (And the induction process stops.)
If $t_k \leq n-N$ then we set $t_{k+1} := t_k + m(y_{i_k}) (<n)$ and choose a point $y_{i_{k+1}}$ with 
$T^{t_{k+1}} x\in V(y_{i_{k+1}})$. This process eventually stops.

From the above construction, for $1\leq k \leq r-2$ 
\begin{equation*}
  \begin{split}
    \#\left(T^{t_k} \left(B_\delta(x, d_n)\right), \, d_{t_{k+1}-t_k}, \frac{\varepsilon}{2}\right)
   &  =  \#\left( T^{t_k}\left(B_\delta(x, d_n)\right), d_{m(y_{i_k})}, \frac{\varepsilon}{2}\right)  \\
   & < e^{(a+\beta)m(y_{i_k})} = e^{(a+\beta)(t_{k+1}-t_k)}. 
   \end{split}
\end{equation*}
For $k=0$ or $k=r-1$, we have $t_{k+1}-t_k \leq N$ and hence 
\[ \#\left(T^{t_k} \left(B_\delta(x, d_n)\right), \, d_{t_{k+1}-t_k}, \, \frac{\varepsilon}{2}\right)
    \leq \#\left(X, d_N, \frac{\varepsilon}{2}\right). \]
Now we use Lemma \ref{lemma: coding each piece} and get 
\begin{equation*}
   \begin{split}
     \#\left(B_\delta(x,d_n), d_n, \varepsilon\right)  & \leq 
     \prod_{k=0}^{r-1} \#\left(T^{t_k}\left(B_\delta(x, d_n)\right), \, d_{t_{k+1}-t_k}, \, \frac{\varepsilon}{2}\right) \\
    & <  \left\{\#\left(X, d_N, \frac{\varepsilon}{2}\right) \right\}^2 \cdot 
       \prod_{k=1}^{r-2} e^{(a+\beta)(t_{k+1}-t_k)} \\
    & \leq  \left\{\#\left(X, d_N, \frac{\varepsilon}{2}\right) \right\}^2 \cdot e^{(a+\beta)n}.
   \end{split}
\end{equation*}
We can regard the term $\left\{\#\left(X, d_N, \frac{\varepsilon}{2}\right) \right\}^2$ as a positive constant
depending on $\varepsilon, \delta, \beta$.
This proves the statement.
\end{proof}

Now we are ready to prove the main result.
We write the statement again.

\begin{theorem}[$=$ Theorem \ref{theorem: main result}]
For any positive number $\delta$
   \begin{equation*}
       \begin{split}
      \overline{\mdim}_M(X, T, d) &= \limsup_{\varepsilon\to 0} 
   \frac{\sup_{x\in X}S\left(B_\delta(x, d_\mathbb{Z}), \varepsilon\right)}{\log(1/\varepsilon)}, \\
     \underline{\mdim}_M(X, T, d) & = \liminf_{\varepsilon\to 0} 
   \frac{\sup_{x\in X}S\left(B_\delta(x,d_{\mathbb{Z}}), \varepsilon\right)}{\log(1/\varepsilon)}.
 \end{split} 
   \end{equation*}
\end{theorem}

\begin{proof}
The proof is simiar to the proof of \cite[Theorem 2.4]{Bowen}.
Let $\varepsilon$ be a positive number.
Set 
\[  a := \sup_{x\in X} S\left(B_\delta(x, d_\mathbb{Z}), \frac{\varepsilon}{4}\right). \]
Let $n$ be a natural number and
$\beta$ a positive number. From Proposition \ref{prop: Bowen ball}
\[  \sup_{x\in X} \#\left(B_\delta(x, d_n), d_n, \varepsilon\right) \leq C \cdot e^{(a+\beta)n}. \]
We choose points $x_1, \dots, x_{M_n}$ in $X$ with 
\[ X = \bigcup_{m=1}^{M_n} B_\delta(x_m, d_n), \quad M_n = \#\left(X, d_n, \delta\right). \]
Then 
\[  \#(X, d_n,\varepsilon) \leq \sum_{m=1}^{M_n} \#\left(B_\delta(x_m, d_n), d_n, \varepsilon\right) 
     \leq C M_n \cdot e^{(a+\beta)n}. \]
Namely 
\[ \log \#(X, d_n,\varepsilon)  \leq \log C + \log \#\left(X, d_n, \delta\right)  + (a+\beta)n. \]
Divide this by $n$ and let $n\to \infty$. Noting that $C$ is independent of $n$, we get 
\[ S(X,\varepsilon) \leq S(X, \delta) + a+\beta. \]
Let $\beta\to 0$:
\[ S(X, \varepsilon) \leq S(X,\delta) + a = S(X, \delta) +  \sup_{x\in X} S\left(B_\delta(x, d_\mathbb{Z}), \frac{\varepsilon}{4}\right). \]
Divide this by $\log (1/\varepsilon)$ and let $\varepsilon\to 0$. We get the statement.
(Notice that $S(X, \delta) \leq \log \#(X, d, \frac{\delta}{2}) < \infty$
by Lemma \ref{lemma: coding each piece}.)
\end{proof}

\section{Generalization to $\mathbb{R}^D$-actions}  \label{section: generalization}

\subsection{Statement of the result}  \label{subsection: statement of the result for R^D}

In the preceding sections, a ``dynamical system" means a compact metrizable space $X$ with 
a homeomorphism $T:X \to X$.
In other words, it is a continuous action of the group $\mathbb{Z}$ on a compact metrizable space.
In order to broaden the applicability of Theorem \ref{theorem: main result}, we would like to generalize it 
to more general group actions. This is essential for studying various geometric examples 
\cite{Gromov, Matsuo--Tsukamoto Brody curves, Tsukamoto Brody}.

Here we consider actions of the group $\mathbb{R}^D$ where $D$ is a natural number.
Probably this is the most basic case\footnote{The argument also works for $\mathbb{Z}^D$-actions.} 
for geometric applications.
$\mathbb{R}^D$ has the standard topology.

A pair $(X, T)$ is called an \textbf{$\mathbb{R}^D$-action} if $X$ is a compact metrizable space and 
$T:\mathbb{R}^D\times X\to X$ is a continuous action.
The (topological) mean dimension of an $\mathbb{R}^D$-action $(X, T)$ is 
denoted as $\mdim(X, T)$. (We do not provide its definition because the topological mean dimension is not 
the main object of this paper.)

Let $(X, T)$ be an $\mathbb{R}^D$-action with a metric $d$ on $X$. 
For a subset $\Omega \subset \mathbb{R}^D$ we define a metric $d_\Omega$ on $X$ by 
\[ d_\Omega(x, y) := \sup_{a\in \Omega} d\left(T^a x, T^a y\right). \]
For $\varepsilon>0$ and $K\subset X$ we define the \textbf{entropy of $K$ at the scale $\varepsilon$}
by 
\[ S(K,\varepsilon) := \limsup_{L\to \infty} \frac{1}{L^D} \log\#\left(K, d_{[0,L]^D},\varepsilon\right). \]
We define the \textbf{upper and lower metric mean dimensions} by 
\[   \overline{\mdim}_M(X, T, d) := \limsup_{\varepsilon\to 0} \frac{S(X, \varepsilon)}{\log(1/\varepsilon)}, \quad 
     \underline{\mdim}_M(X, T, d)  := \liminf_{\varepsilon\to 0} \frac{S(X, \varepsilon)}{\log(1/\varepsilon)}. \]
They bound the mean dimension from above \cite[Theorem 4.2]{Lindenstrauss--Weiss}:
\begin{equation} \label{eq: metric mean dimension bound mean dimension for R^D actions}
    \mdim(X, T) \leq  \underline{\mdim}_M(X, T, d)  \leq \overline{\mdim}_M(X, T, d). 
\end{equation}

For $x\in X$, $\delta>0$ and a subset $\Omega\subset \mathbb{R}^D$, we set 
\[  B_\delta\left(x, d_\Omega\right) := \{y\in X|\, d_\Omega(x, y) \leq \delta\}. \]
In particular 
\[ B_\delta\left(x, d_{\mathbb{R}^D}\right) = \{y\in X|\, d_{\mathbb{R}^D}(x, y) \leq \delta\}
     = \{y\in X|\, \forall a\in \mathbb{R}^D: \> d(T^a x, T^a y) \leq \delta\}. \]
The following is the main result of this section:

\begin{theorem}  \label{theorem: main result for R^D}
Let $(X, T)$ be an $\mathbb{R}^D$-action with a metric $d$.
For any $\delta>0$
   \begin{equation*}
       \begin{split}
      \overline{\mdim}_M(X, T, d) &= \limsup_{\varepsilon\to 0} 
   \frac{\sup_{x\in X}S\left(B_\delta(x, d_{\mathbb{R}^D}), \varepsilon\right)}{\log(1/\varepsilon)}, \\
     \underline{\mdim}_M(X, T, d) & = \liminf_{\varepsilon\to 0} 
   \frac{\sup_{x\in X}S\left(B_\delta(x,d_{\mathbb{R}^D}), \varepsilon\right)}{\log(1/\varepsilon)}.
 \end{split} 
   \end{equation*}
\end{theorem}

We would like to mention a related result.
In the study of the mean dimension of the system of Brody curves (see \S \ref{section: example}) in \cite{Tsukamoto Brody},
the author encountered a problem of how to formulate the local nature of metric mean dimension.
At that time the author did not realize the paper of Bowen \cite{Bowen}.
In \cite[Lemma 2.5]{Tsukamoto Brody}, the following lemma was proved.

\begin{lemma}[\cite{Tsukamoto Brody}] \label{lemma: weaker version}
Let $(X, T)$ be an $\mathbb{R}^D$-action with a metric $d$. For any $\delta>0$ and $R\geq 0$, 
\[  \overline{\mdim}_M\left(X, T, d\right) = \limsup_{\varepsilon \to 0} \left(\limsup_{L\to \infty}
    \frac{\sup_{x\in X}\log\#\left(B_\delta\left(x, d_{[-R, L+R]^D}\right), d_{[0,L]^D}, \varepsilon\right)}{L^D \cdot \log(1/\varepsilon)}\right). \]
We also have a similar result for the lower metric mean dimension.
\end{lemma}

This is an important technical ingredient of \cite{Tsukamoto Brody}. 
The basic philosophy behind this lemma is the same as in Theorem \ref{theorem: main result for R^D}.
It utilizes the local nature of metric mean dimension. We can reduce the global problem to a more local problem of 
studying the Bowen balls $B_\delta\left(x, d_{[-R, L+R]^D}\right)$.
This was very useful for studying the mean dimension of the system of Brody curves.

Lemma \ref{lemma: weaker version} easily follows from Theorem \ref{theorem: main result for R^D}.
So Theorem \ref{theorem: main result for R^D} is a stronger result.
Moreover the statement of Theorem \ref{theorem: main result for R^D} is substantially simpler than 
Lemma \ref{lemma: weaker version}.
So the author thinks that it provides a ``right'' formulation of the local nature of metric mean dimension.

\subsection{Tiling argument}  \label{subsection: covering argument}

Conceptually, the proof of Theorem \ref{theorem: main result for R^D} is the same as the 
proof of Theorem \ref{theorem: main result}. 
But technically it is a bit more complicated. 
We need a kind of \textit{tiling argument}
originally due to Ornstein--Weiss \cite{Ornstein--Weiss}.
This subsection is a preparation for it.

We consider the $\ell^\infty$-norm on $\mathbb{R}^D$:
\[ \left|(x_1, \dots, x_D)\right|_\infty := \max_{1\leq n \leq D} |x_n|. \]
We always think that $\mathbb{R}^D$ is endowed with this norm (not the $\ell^2$-norm).
In this section a ``cube" $\Lambda$ in $\mathbb{R}^D$ means a set of the form 
\[ \Lambda = u + [0,L]^D = \{u+x|\, x\in [0,L]^D\}, \]
where $u\in \mathbb{R}^D$ and $L>0$.
We set $\ell(\Lambda) := L$.
We also set 
\[  3\Lambda := u + [-L, 2L]^D. \]  
We denote by $\vol(\cdot)$ the $D$-dimensional Lebesgue measure on $\mathbb{R}^D$.

The following is a kind of (finite) Vitali covering lemma.

\begin{lemma} \label{lemma: Vitali}
Let $\Lambda_1, \dots, \Lambda_N$ be cubes in $\mathbb{R}^D$. 
There is a disjoint subfamily $\{\Lambda_{n_1}, \dots, \Lambda_{n_k}\}$ such that 
\[  \Lambda_1\cup \Lambda_2 \cup \dots \cup \Lambda_N \subset 
     3\Lambda_{n_1} \cup 3\Lambda_{n_2} \cup \dots \cup 3\Lambda_{n_k}. \]
In particular, 
\[ \vol\left(\Lambda_{n_1}\cup \dots \cup \Lambda_{n_k}\right) \geq 
    3^{-D} \vol\left(\Lambda_1\cup \dots \cup \Lambda_N\right). \]
\end{lemma}

\begin{proof}
We use a greedy algorithm.
Let $\Lambda_{n_1}$ be one of the cubes in $\{\Lambda_1, \dots, \Lambda_N\}$ which have the maximum side length $\ell(\Lambda_n)$.
We discard all the cubes which have non-empty intersection with $\Lambda_{n_1}$.
Next let $\Lambda_{n_2}$ be one of the remaining cubes which have the maximum side length.  
Then we discard all the cubes which have non-empty intersection with $\Lambda_{n_2}$.
We continue this process. It eventually stops and we get a disjoint family $\{\Lambda_{n_1}, \dots, \Lambda_{n_k}\}$.

Let $\Lambda_m$ be an arbitrary cube in the initial family. If it is not chosen in the process, then there is $\Lambda_{n_i}$ satisfying 
\[  \Lambda_m \cap \Lambda_{n_i} \neq \emptyset, \quad \ell(\Lambda_m) \leq \ell(\Lambda_{n_i}). \]
Then $\Lambda_m\subset 3\Lambda_{n_i}$. 
\end{proof}

For $r>0$ and a subset $\Omega \subset \mathbb{R}^D$, we define the \textbf{$r$-interior} and \textbf{$r$-boundary} of 
$\Omega$ by 
\[  \mathrm{int} (\Omega,r) := \{x\in \Omega|\, x+[-r,r]^D\subset \Omega\}, \]
\[    \partial(\Omega, r) := \{x\in \mathbb{R}^D|\, \exists y \in \Omega, \, \exists z\in \mathbb{R}^D\setminus \Omega: 
      |x-y|_\infty\leq r, \, |x-z|_\infty \leq r \}. \] 
We also set 
\[ B_r(\Omega) := \Omega \cup \partial(\Omega, r) = \{x\in \mathbb{R}^D|\, \exists y\in \Omega: \> |x-y|_\infty \leq r\}
     = \mathrm{int}(\Omega, r) \cup \partial(\Omega, r). \]
For $\Omega, \Omega'\subset \mathbb{R}^D$, we have
\[  \partial(\Omega\cup \Omega', r) \subset \partial(\Omega, r) \cup \partial(\Omega', r). \]

Let $\mathcal{C} = \{\Lambda_1, \dots, \Lambda_N\}$ be a finite family of cubes in $\mathbb{R}^D$.
We set 
\[ \ell_{\max}(\mathcal{C}) := \max_{1\leq n\leq N} \ell(\Lambda_n), \quad 
    \ell_{\min}(\mathcal{C}) := \min_{1\leq n\leq N} \ell(\Lambda_n). \]

The next proposition is a main technical ingredient of the proof of Theorem \ref{theorem: main result for R^D}.

\begin{proposition}  \label{prop: covering lemma}
For any $\eta>0$ there exists a natural number $K = K(\eta)>0$ satisfying the following statement.
Let $\Omega \subset \mathbb{R}^D$ be a bounded Borel subset, and let $\mathcal{C}_1, \dots, \mathcal{C}_K$ be 
finite families of cubes in $\mathbb{R}^D$. We assume:
\begin{enumerate}
   \item $\ell_{\max}(\mathcal{C}_1) \geq 1$ and $\ell_{\min}(\mathcal{C}_{k+1}) \geq K \cdot \ell_{\max}(\mathcal{C}_k)$ for all 
            $1\leq k\leq K-1$. 
   \item $\vol\left\{\partial\left(\Omega, \ell_{\max}(\mathcal{C}_K)\right)\right\} < \frac{\eta}{3} \cdot \vol(\Omega)$.
   \item $\Omega \subset \bigcup_{\Lambda\in \mathcal{C}_k} \Lambda$ 
            for every $1\leq k\leq K$.
\end{enumerate}
Then there exists a disjoint subfamily $\mathcal{A}\subset \mathcal{C}_1\cup \dots \cup \mathcal{C}_K$ satisfying 
\[ \bigcup_{\Lambda\in \mathcal{A}}\Lambda  \subset \Omega, \quad 
 \vol\left\{B_1\left(\Omega \setminus \bigcup_{\Lambda\in \mathcal{A}}\Lambda\right) \right\} < \eta \cdot \vol(\Omega). \]
\end{proposition}

This is a rather complicated statement. Before proving the lemma, we explain its meaning more intuitively.
First, since $K$ will be assumed very large, the number $K$ of the given families $\mathcal{C}_1, \dots, \mathcal{C}_K$ is very large.
The condition (1) means that that every cube in $\mathcal{C}_{k+1}$ is 
much larger than cubes in $\mathcal{C}_k$.
(So we have $K$ different scales.)
The condition (2) means that the set $\Omega$ is much larger than cubes in $\mathcal{C}_1\cup \dots \cup \mathcal{C}_K$.
The condition (3) is simple.
The conclusion means that a large portion of $\Omega$ can be covered by mutually disjoint cubes in 
$\mathcal{C}_1\cup \dots \cup \mathcal{C}_K$.

\begin{proof}
We assume that $K>1$ is sufficiently large so that 
\begin{itemize}
   \item $K \cdot 3^{-D}\cdot \frac{\eta}{3}>1$.
   \item For $r>0$ and a cube $\Lambda\subset \mathbb{R}^D$, if $\ell(\Lambda) \geq K r$ then 
            $\vol\left(\partial(\Lambda, r)\right) < \frac{\eta}{3}\cdot \vol(\Lambda)$.
\end{itemize}

We will inductively construct $\mathcal{A}_k \subset \mathcal{C}_{K+1-k}$ for 
$k=1, 2, \dots, K$ and finally define 
$\mathcal{A} = \mathcal{A}_1\cup \dots \cup \mathcal{A}_K$.
We also set $\Omega_k = 
\Omega \setminus \bigcup_{\Lambda\in \mathcal{A}_1\cup \dots \cup \mathcal{A}_k} \Lambda$.
(Let $\Omega_0 := \Omega$.)

Suppose we have defined $\mathcal{A}_1, \dots, \mathcal{A}_k$. We are going to define $\mathcal{A}_{k+1}$.
(The initial step of the induction is the same; set $k=0$ in the argument below.)
Set 
\[ \mathcal{C}'_{K-k} = \left\{\Lambda\in \mathcal{C}_{K-k}\middle|\, 
    \Lambda\cap \mathrm{int}\left(\Omega_k,\ell_{\max}(\mathcal{C}_{K-k})\right) \neq \emptyset\right\}. \]
We have 
\[ \mathrm{int}\left(\Omega_k, \ell_{\max}(\mathcal{C}_{K-k})\right)
  \subset \bigcup_{\Lambda\in \mathcal{C}'_{K-k}} \Lambda \subset \Omega_k. \]
In particular, every cube in $\mathcal{C}'_{K-k}$ has no intersection with 
cubes in $\mathcal{A}_1\cup \dots \cup \mathcal{A}_k$. 

By Lemma \ref{lemma: Vitali} there exists a disjoint subfamily $\mathcal{A}_{k+1}\subset \mathcal{C}'_{K-k}$
satisfying 
\[ \vol\left(\bigcup_{\Lambda\in \mathcal{A}_{k+1}} \Lambda\right) 
    \geq 3^{-D} \cdot \vol\left(\bigcup_{\Lambda\in \mathcal{C}'_{K-k}} \Lambda\right) 
    \geq 3^{-D} \cdot \vol\left\{\mathrm{int}\left(\Omega_k, \ell_{\max}(\mathcal{C}_{K-k})\right)\right\}. \]

Now we have defined $\mathcal{A}_k$ for all $1\leq k\leq K$. 
We set $\mathcal{A} = \mathcal{A}_1\cup \dots \cup \mathcal{A}_K$.
This is a disjoint family of cubes and $\bigcup_{\Lambda\in \mathcal{A}}\Lambda  \subset \Omega$.

\begin{claim}
There exists $0\leq k <K$ satisfying 
\begin{equation} \label{eq: choice of k in the covering lemma}
    \vol\left\{\mathrm{int}\left(\Omega_k, \ell_{\max}(\mathcal{C}_{K-k})\right)\right\} < \frac{\eta}{3} \cdot \vol(\Omega). 
\end{equation}
\end{claim}

\begin{proof}
Suppose this is false. Then for all $0\leq k < K$
\[  \vol\left(\bigcup_{\Lambda\in \mathcal{A}_{k+1}} \Lambda\right) \geq 
     3^{-D}\cdot \vol\left\{\mathrm{int}\left(\Omega_k, \ell_{\max}(\mathcal{C}_{K-k})\right)\right\} 
    \geq  3^{-D}\cdot \frac{\eta}{3} \cdot \vol(\Omega). \]
Then 
\[ \vol\left(\Omega_K\right) \leq \vol(\Omega) - K \cdot 3^{-D} \cdot \frac{\eta}{3} \cdot \vol(\Omega)
    = \left(1 - K\cdot 3^{-D}\cdot \frac{\eta}{3}\right) \vol(\Omega) < 0 \]
because we assumed $K \cdot 3^{-D}\cdot \frac{\eta}{3}>1$. This is a contradiction.
\end{proof}

Let $0\leq k <K$ be the integer satisfying (\ref{eq:  choice of k in the covering lemma}).
Set $r = \ell_{\max}(\mathcal{C}_{K-k}) \geq 1$.
We have 
\[ \partial(\Omega_k, r) \subset  \partial(\Omega, r) \cup 
    \bigcup_{\Lambda\in \mathcal{A}_1\cup \dots \cup \mathcal{A}_{k}} \partial(\Lambda, r). \]
From the choice of $K$ in the beginning, for $\Lambda\in \mathcal{A}_1\cup \dots \cup \mathcal{A}_{k}$
\[  \vol\left(\partial(\Lambda, r)\right) < \frac{\eta}{3} \cdot \vol(\Lambda). \]
Then 
\[ \vol\left(\partial(\Omega_k, r)\right) \leq \vol\left(\partial (\Omega, r)\right) + 
     \sum_{\Lambda\in \mathcal{A}_1\cup \dots \cup \mathcal{A}_k} \frac{\eta}{3}\cdot \vol(\Lambda) 
     < \frac{\eta}{3}\vol(\Omega) + \frac{\eta}{3} \vol(\Omega) = \frac{2\eta}{3} \cdot \vol(\Omega). \]
Then we estimate the volume of $B_r(\Omega_k) = \mathrm{int}(\Omega_k, r) \cup \partial(\Omega_k, r)$ by
\[  \vol\left(B_r(\Omega_k)\right) < \frac{\eta}{3}\cdot \vol(\Omega) + \frac{2\eta}{3}\cdot \vol(\Omega) 
      = \eta\cdot \vol(\Omega). \]
Since $B_1(\Omega_K) \subset B_r(\Omega_k)$, this proves the statement.
\end{proof}

\subsection{Proof of Theorem \ref{theorem: main result for R^D}}  \label{subsection: proof of main result for R^D}

Throughout this subsection we assume that $(X, T)$ is an $\mathbb{R}^D$-action with a metric $d$.

\begin{lemma} \label{lemma: coding small pieces}
Let $\varepsilon>0$ and $F\subset X$. Let $\Omega, \Omega_1, \dots, \Omega_N\subset \mathbb{R}^D$ be bounded subsets.
If $\Omega \subset \Omega_1\cup \dots \cup \Omega_N$ then 
\[  \#(F, d_\Omega, \varepsilon) \leq 
     \prod_{n=1}^N \#\left(F, d_{\Omega_n}, \frac{\varepsilon}{2}\right). \]
\end{lemma}

\begin{proof}
Let $E_n \subset X$ be a $\left(d_{\Omega_n}, \frac{\varepsilon}{2}\right)$-spanning set of $F$.
For each $(x_1, \dots, x_N)$ in $E_1\times \dots \times E_N$, we set 
\[  V(x_1, \dots, x_N) := \left\{y\in X\middle|\, \forall 1\leq n\leq N: \, d_{\Omega_n}(y, x_n) \leq \frac{\varepsilon}{2} \right\}. \]
The diameter of this set with respect to $d_\Omega$ is less than or equal to $\varepsilon$.
The sets $V(x_1, \dots, x_N)$, $(x_1, \dots, x_N)\in E_1\times \dots \times E_N$, cover $F$.
We pick a point from each (non-empty) $V(x_1, \dots, x_N)$. Then we get a $(d_\Omega, \varepsilon)$-spanning set of $F$.
\end{proof}

\begin{lemma}  \label{lemma: trivial bound}
For any $\varepsilon>0$ and any bounded Borel subset $\Omega\subset \mathbb{R}^D$
\[  \#\left(X, d_\Omega,\varepsilon\right) \leq \left(\#\left(X, d_{[0,1]^D},\frac{\varepsilon}{2}\right)\right)^{\vol \left(B_1(\Omega)\right)}. \]
\end{lemma}

\begin{proof}
Let $A$ be the set of $u\in \mathbb{Z}^D$ satisfying $\left(u+[0,1]^D\right)\cap \Omega \neq \emptyset$. Then 
\[ \Omega \subset \bigcup_{u\in A} \left(u+[0,1]^D\right) \subset B_1(\Omega). \]
Then $|A| \leq \vol\left(B_1(\Omega)\right)$ and (by Lemma \ref{lemma: coding small pieces}) 
\[  \#\left(X, d_\Omega, \varepsilon\right) \leq \prod_{u\in A} \#\left(X, d_{u + [0,1]^D}, \frac{\varepsilon}{2}\right). \]
We have 
\[  \#\left(X, d_{u + [0,1]^D}, \frac{\varepsilon}{2}\right) = \#\left(X, d_{[0,1]^D}, \frac{\varepsilon}{2}\right). \]
This proves the statement.
\end{proof}

\begin{proposition} \label{prop: key for main result for R^D}
Let $\varepsilon, \delta, \beta$ be positive numbers.
There exists a positive number $N=N(\varepsilon, \delta, \beta)>0$ satisfying the following statement.
Set
\[  a := \sup_{x\in X} S\left(B_\delta(x, d_{\mathbb{R}^D}), \frac{\varepsilon}{4}\right). \]
Then 
\[  \sup_{x\in X} \#\left(B_\delta\left(x, d_{[-N, L+N]^D}\right), d_{[0,L]^D}, \varepsilon\right) 
     < e^{(a+\beta)L^D} \]
for all sufficiently large $L>1$.
\end{proposition}

\begin{proof}
We choose $\eta>0$ so small that 
\begin{equation} \label{eq: choice of eta} 
  \left(\#\left(X, d_{[0,1]^D},\frac{\varepsilon}{4}\right)\right)^\eta < e^{\frac{\beta}{2}}. 
\end{equation}
Let $K=K(\eta)>1$ be the positive number for this $\eta$ introduced in Proposition \ref{prop: covering lemma}.

For $k=1, 2, \dots, K$, we will inductively construct a finite subset $Y_k\subset X$ and, for each $y\in Y_k$,
positive numbers $L_k(y), M_k(y)$ and open neighborhoods $U_k(y), V_k(y)$ of $y$ satisfying the following.
\begin{itemize}
  \item $L_1(y) >1$ and $L_k(y) >K \max_{z\in Y_{k-1}} L_{k-1}(z)$ for $k\geq 2$.
  \item $\#\left(U_k(y), d_{[0, L_k(y)]^D}, \frac{\varepsilon}{2}\right) < \exp\left(\left(a+\frac{\beta}{2}\right)L_k(y)^D\right)$. 
  \item For every $v\in V_k(y)$, $B_\delta\left(v, d_{[-M_k(y), M_k(y)]^D}\right) \subset U_k(y)$.
  \item $X = \bigcup_{y\in Y_k} V_k(y)$ for every $1\leq k \leq K$.
\end{itemize}

Suppose we have defined the data for the $(k-1)$-th step. 
We are going to construct the data for the $k$-th step. 
(The initial step, $k=1$, can be treated in the same way.)
Let $y\in X$ be an arbitrary point. Since we have 
\[  S\left(B_\delta\left(y, d_{\mathbb{R}^D}\right), \frac{\varepsilon}{4}\right) \leq a < a+\frac{\beta}{2}, \]
we can choose $L_k(y) > K \max_{z\in Y_k} L_{k-1}(z)$ (when $k=1$, we assume 
$L_1(y) >1$) satisfying 
\[ \frac{1}{L_k(y)^D} \log \#\left(B_\delta\left(y, d_{\mathbb{R}^D}\right), d_{[0,L_k(y)]^D}, \frac{\varepsilon}{4}\right)
     < a + \frac{\beta}{2}. \]
Then there is $E_k(y)\subset X$ which is a $\left(d_{[0, L_k(y)]^D}, \frac{\varepsilon}{4}\right)$-spanning 
set of $B_\delta\left(y, d_{\mathbb{R}^D}\right)$ with 
\[  |E_k(y)| < \exp\left(\left(a+\frac{\beta}{2}\right) |L_k(y)|^D\right). \]
We define
\[ U_k(y) := \left\{x\in X\middle|\, \exists z\in E_k(y): d_{[0, L_k(y)]^D}(x, z) < \frac{\varepsilon}{2} \right\}. \]
This is an open set containing $B_\delta\left(y, d_{\mathbb{R}^D}\right)$ with 
\[ \#\left(U_k(y), d_{[0, L_k(y)]^D}, \frac{\varepsilon}{2}\right) \leq |E_k(y)|
     < \exp\left(\left(a+\frac{\beta}{2}\right) |L_k(y)|^D\right). \]
Since 
\[  B_\delta\left(y, d_{\mathbb{R}^D}\right) = \bigcap_{M=1}^\infty B_\delta\left(y, d_{[-M, M]^D}\right), \]
there exists $M_k(y)>0$ satisfying 
\[ B_\delta\left(y, d_{[-M_k(y), M_k(y)]^D}\right)  \subset U_k(y). \]
If $\delta'>\delta$ is sufficiently close to $\delta$ then we also have 
\[  B_{\delta'}\left(y, d_{[-M_k(y), M_k(y)]^D}\right)  \subset U_k(y). \]
Set 
\[ V_k(y) := \left\{x\in X\middle|\, d_{[-M_k(y), M_k(y)]^D}(x, y) < \delta'-\delta\right\}. \]
Then for every $v\in V_k(y)$ we have $B_\delta\left(v, d_{[-M_k(y), M_k(y)]^D}\right)  \subset U_k(y)$.
$V_k(y)$ is an open neighborhood of $y$. So we choose a finite set $Y_k\subset X$ such that 
$V_k(y)$, $y\in Y_k$, cover $X$.
We have finished the construction of the $k$-th step.
So the induction works.

We fix 
\begin{equation}  \label{eq: choice of N in the main proposition for R^D}
   N > \max\{M_k(y)|\, 1\leq k\leq K, y\in Y_k\}. 
\end{equation}
We assume that $L>1$ is any sufficiently large number such that the cube $\Omega :=[0, L]^D$ satisfies 
\[ \vol\left(\partial\left(\Omega, \max_{y\in Y_K} L_K(y)\right)\right) < \frac{\eta}{3} \cdot \vol(\Omega). \]
We are going to prove 
\[ \sup_{x\in X} \#\left(B_\delta\left(x, d_{[-N,L+N]^D}\right), d_{[0,L]^D}\right) < e^{(a+\beta)L^D}. \]

Take any $x\in X$. 
For each $1\leq k\leq K$ and $t\in \mathbb{Z}^D\cap \Omega$, we pick
$y\in Y_k$ with $T^t x\in V_k(y)$.
Set $\Lambda_{k, t} := t + [0, L_k(y)]^D$.
From the choice of $N$ in (\ref{eq: choice of N in the main proposition for R^D}), 
\[ T^t\left(B_\delta\left(x, d_{[-N, L+N]^D}\right)\right) \subset 
    B_\delta\left(T^t x, d_{[-M_k(y), M_k(y)]^D}\right) \subset U_k(y). \]
Hence 
\begin{equation}  \label{eq: evaluate spanning number over cubes}
   \#\left(B_\delta\left(x, d_{[-N, L+N]^D}\right), d_{\Lambda_{k, t}}, \frac{\varepsilon}{2}\right)
   \leq \#\left(U_k(y), d_{[0, L_k(y)]^D}, \frac{\varepsilon}{2}\right) 
   < e^{\left(a+\frac{\beta}{2}\right) \vol(\Lambda_{k, t})}. 
\end{equation}

Let $\mathcal{C}_k$ be the set of cubes $\Lambda_{k, t}$ $(t\in \mathbb{Z}^D\cap \Omega)$.
(Notice that this depends on $x\in X$.) 
The cubes in $\mathcal{C}_k$ cover $\Omega = [0, L]^D$.

Now we apply Proposition \ref{prop: covering lemma} to $\Omega$ and $\mathcal{C}_1, \dots, \mathcal{C}_K$.
Then there is a disjoint subfamily $\mathcal{A}\subset \mathcal{C}_1\cup \dots \cup \mathcal{C}_K$
such that (set $\Omega' := \Omega \setminus \bigcup_{\Lambda\in \mathcal{A}} \Lambda$)
\[  \bigcup_{\Lambda\in \mathcal{A}} \Lambda \subset \Omega, \quad 
    \vol\left( B_1\left(\Omega'\right)\right) < \eta \cdot \vol(\Omega). \]
From Lemma \ref{lemma: trivial bound} and the choice of $\eta$ in (\ref{eq: choice of eta})
\begin{equation*}
  \begin{split}
    \#\left(B_\delta\left(x, d_{[-N, L+N]^D}\right), d_{\Omega'}, \frac{\varepsilon}{2}\right)
     & \leq   \#\left(X, d_{\Omega'}, \frac{\varepsilon}{2}\right)    \\
     &  \leq  \left(\#\left(X, d_{[0,1]^D},\frac{\varepsilon}{4}\right)\right)^{\vol\left( B_1\left(\Omega'\right)\right)} \\
     & \leq  \left(\#\left(X, d_{[0,1]^D},\frac{\varepsilon}{4}\right)\right)^{\eta \cdot \vol(\Omega)} \\
     & < e^{\frac{\beta}{2} \vol(\Omega)}. 
   \end{split}
\end{equation*}
On the other hand, for each $\Lambda\in \mathcal{A}$, from (\ref{eq: evaluate spanning number over cubes})
\[   \#\left(B_\delta\left(x, d_{[-N, L+N]^D}\right), d_{\Lambda}, \frac{\varepsilon}{2}\right)
    < e^{\left(a+\frac{\beta}{2}\right) \vol(\Lambda)}. \]
From Lemma \ref{lemma: coding small pieces}
\begin{equation*}
  \begin{split}
   & \#\left(B_\delta\left(x, d_{[-N, L+N]^D}\right), d_{[0,L]^D}, \varepsilon\right) \\
   & \leq 
   \#\left(B_\delta\left(x, d_{[-N, L+N]^D}\right), d_{\Omega'}, \frac{\varepsilon}{2}\right) \cdot 
    \prod_{\Lambda\in \mathcal{A}} \#\left(B_\delta\left(x, d_{[-N, L+N]^D}\right), d_{\Lambda}, \frac{\varepsilon}{2}\right) \\
   & < e^{\frac{\beta}{2} \vol(\Omega)} \cdot \prod_{\Lambda\in \mathcal{A}} e^{\left(a+\frac{\beta}{2}\right) \vol(\Lambda)} \\
   & \leq e^{(a+\beta) \vol(\Omega)}  = e^{(a+\beta)L^D}.
  \end{split}
\end{equation*}
This holds for any $x\in X$. So we have proved the statement of the proposition.
\end{proof}

Now we are ready to prove Theorem \ref{theorem: main result for R^D}.
We write the statement again.

\begin{theorem}[$=$ Theorem \ref{theorem: main result for R^D}]
For any $\delta>0$
   \begin{equation*}
       \begin{split}
      \overline{\mdim}_M(X, T, d) &= \limsup_{\varepsilon\to 0} 
   \frac{\sup_{x\in X}S\left(B_\delta(x, d_{\mathbb{R}^D}), \varepsilon\right)}{\log(1/\varepsilon)}, \\
     \underline{\mdim}_M(X, T, d) & = \liminf_{\varepsilon\to 0} 
   \frac{\sup_{x\in X}S\left(B_\delta(x,d_{\mathbb{R}^D}), \varepsilon\right)}{\log(1/\varepsilon)}.
 \end{split} 
   \end{equation*}
\end{theorem}

\begin{proof}
Let $\varepsilon$ and $\beta$ be any positive numbers.
Let $N=N(\varepsilon, \delta, \beta)$ be the positive number given in Proposition 
\ref{prop: key for main result for R^D}. Set 
\[ a:= \sup_{x\in X} S\left(B_\delta\left(x, d_{\mathbb{R}^D}\right), \frac{\varepsilon}{4}\right). \]
For any $L>0$, there are $x_1, \dots, x_M\in X$ such that 
\[ X = \bigcup_{m=1}^M B_\delta\left(x_m, d_{[-N,L+N]^D}\right), \quad 
    M = \#\left(X, d_{[-N, L+N]^D}, \delta\right) = \#\left(X, d_{[0, L+2N]^D}, \delta\right). \]
Then
\begin{equation*}
   \begin{split}
    \#\left(X, d_{[0, L]^D}, \varepsilon\right) & \leq \sum_{m=1}^M 
    \#\left(B_\delta\left(x_m, d_{[-N,L+N]^D}\right), d_{[0, L]^D}, \varepsilon\right)  \\
    & \leq M \cdot \sup_{x\in X}   \#\left(B_\delta\left(x, d_{[-N,L+N]^D}\right), d_{[0, L]^D}, \varepsilon\right) \\
    & < M \cdot e^{(a+\beta)L^D}
   \end{split}
\end{equation*}
for all sufficiently large $L$ by Proposition \ref{prop: key for main result for R^D}.
So, for large $L$
\[ \log \#\left(X, d_{[0, L]^D}, \varepsilon\right) \leq \log \#\left(X, d_{[0, L+2N]^D}, \delta\right)
    + (a+\beta)L^D. \]
Divide this by $L^D$ and let $L\to \infty$.
\[  S(X, \varepsilon) \leq S(X, \delta) + a + \beta. \]
Let $\beta\to 0$:
\[ S(X, \varepsilon) \leq S(X,\delta) + a = 
S(X,\delta) + \sup_{x\in X} S\left(B_\delta\left(x, d_{\mathbb{R}^D}\right), \frac{\varepsilon}{4}\right). \]
Divide this by $\log(1/\varepsilon)$ and let $\varepsilon\to 0$. 
Notice that $S(X,\delta) \leq \log\#\left(X,d_{[0,1]^D},\frac{\delta}{2}\right)<\infty$ 
by Lemma \ref{lemma: trivial bound}.
We get the conclusion.
\end{proof}

\section{Example: Brody curves}  \label{section: example}

In this section we explain how to use Theorem \ref{theorem: main result for R^D}
by an example.
Our example is the $\mathbb{C}$-action on the space of \textbf{Brody curves}
(Lipschitz holomorphic curves).
We revisit a result proved in \cite{Tsukamoto Brody}.

\subsection{$\mathbb{C}$-action on the space of Brody curves} 
\label{subsection: C-action on the space of Brody curves}

Let $z=x+y\sqrt{-1}$ be the standard coordinate of the complex plane $\mathbb{C}$.
We consider the complex projective space $\mathbb{C}P^N$ with the Fubini--Study metric.
For a holomorphic map $f:\mathbb{C}\to \mathbb{C}P^N$, we denote by $|df|(z)\geq 0$ the local 
Lipschitz constant at $z\in \mathbb{C}$. Explicitly, for $f=[f_0:f_1:\dots:f_N]$
\[  |df|^2(z) = \frac{1}{4\pi}\left(\frac{\partial^2}{\partial x^2} + \frac{\partial^2}{\partial y^2}\right)
     \log\left(|f_0|^2 + |f_1|^2+\dots+|f_N|^2\right). \]
For $\lambda>0$ we define $\mathcal{M}_\lambda(\mathbb{C}P^N)$ as the space of 
holomorphic maps $f:\mathbb{C}\to \mathbb{C}P^N$ satisfying $|df|\leq \lambda$ 
all over the plane\footnote{A Lipschitz holomorphic map from $\mathbb{C}$ is called a Brody curve.
Brody \cite{Brody} found its importance in the study of Kobayashi hyperbolicity.}.
This is endowed with the compact-open topology and becomes a compact metrizable space.
The group $\mathbb{C} = \mathbb{R}^2$ continuously acts on it by 
\[  T:\mathbb{C}\times \mathcal{M}_\lambda(\mathbb{C}P^N) \to \mathcal{M}_\lambda(\mathbb{C}P^N), \quad 
      \left(a, f(z)\right) \mapsto f(z+a). \]
We would like to study the mean dimension of $\left(\mathcal{M}(\mathbb{C}P^N), T\right)$.
It is easy to see that 
\[ \mdim\left(\mathcal{M}_\lambda(\mathbb{C}P^N), T\right) 
     = \lambda^2\cdot \mdim\left(\mathcal{M}_1(\mathbb{C}P^N), T\right). \]
So it is enough to study the case of $\lambda=1$.
(Nevertheless, it is useful to consider other $\mathcal{M}_\lambda(\mathbb{C}P^N)$
even for the study of $\mathcal{M}_1(\mathbb{C}P^N)$.
See Proposition \ref{prop: resolution} below.)

For $f\in \mathcal{M}_\lambda(\mathbb{C}P^N)$, we define its \textbf{energy density} by 
\[ \rho(f) := \lim_{L\to \infty} \left(\frac{1}{L^2}\sup_{a\in \mathbb{C}} \int_{a+[0, L]^2} |df|^2\, dxdy\right). \]
We define $\rho_\lambda(\mathbb{C}P^N)$ as the supremum of $\rho(f)$ over 
$f\in \mathcal{M}_\lambda(\mathbb{C}P^N)$. We have 
\[  \rho_\lambda(\mathbb{C}P^N)  = \lambda^2 \cdot \rho_1(\mathbb{C}P^N). \]

In \cite{Matsuo--Tsukamoto Brody curves} we proved the lower bound 
\[ \mdim\left(\mathcal{M}_1(\mathbb{C}P^N), T\right)  \geq 2(N+1) \rho_1(\mathbb{C}P^N). \]
On the other hand, it was proved in \cite{Tsukamoto Brody} that
\[ \mdim\left(\mathcal{M}_1(\mathbb{C}P^N), T\right)  \leq 2(N+1) \rho_1(\mathbb{C}P^N). \]
So we get the formula 
\[ \mdim\left(\mathcal{M}_1(\mathbb{C}P^N), T\right)  = 2(N+1) \rho_1(\mathbb{C}P^N). \]

The purpose here is to explain the proof of the upper bound 
$\mdim\left(\mathcal{M}_1(\mathbb{C}P^N), T\right)  \leq 2(N+1) \rho_1(\mathbb{C}P^N)$.
The proof uses metric mean dimension.

For $f, g\in \mathcal{M}_\lambda(\mathbb{C}P^N)$ we define 
\[ d(f, g) := \max_{z\in [0,1]^2} \mathbf{d}_{\mathrm{FS}}\left(f(z), g(z)\right). \]
Here $\mathbf{d}_{\mathrm{FS}}(\cdot, \cdot)$ is the Fubini--Study metric.
This $d(f, g)$ becomes a metric on $\mathcal{M}_\lambda(\mathbb{C}P^N)$ by the unique continuation 
principle. Notice that we have
\[ d_{\mathbb{C}}(f, g) = \sup_{z\in \mathbb{C}} \mathbf{d}_{\mathrm{FS}}\left(f(z), g(z)\right). \]
We would like to prove the upper bound on the upper metric mean dimension
\begin{equation} \label{eq: upper bound on metric mean dimension}
    \overline{\mdim}_M\left(\mathcal{M}_1(\mathbb{C}P^N), T, d\right) \leq 2(N+1) \rho_1(\mathbb{C}P^N). 
\end{equation}
Since metric mean dimension bounds mean dimension (see (\ref{eq: metric mean dimension bound mean dimension for R^D actions})),
it follows from (\ref{eq: upper bound on metric mean dimension}) that 
\[ \mdim\left(\mathcal{M}_1(\mathbb{C}P^N), T\right)  \leq 2(N+1) \rho_1(\mathbb{C}P^N). \]
So the problem is how to prove (\ref{eq: upper bound on metric mean dimension}).
The proof will be given in \S \ref{subsection: proof of the upper bound}.

In \cite{Tsukamoto Brody} the upper bound (\ref{eq: upper bound on metric mean dimension}) was proved 
by using Lemma \ref{lemma: weaker version} (a weaker version of Theorem \ref{theorem: main result for R^D}).
Here we will prove (\ref{eq: upper bound on metric mean dimension}) by using Theorem 
\ref{theorem: main result for R^D}.
The basic structures of the two proofs are the same. But the use of Theorem \ref{theorem: main result for R^D}
makes the argument clean.

By Theorem \ref{theorem: main result for R^D} 
\[  \overline{\mdim}_M\left(\mathcal{M}_1(\mathbb{C}P^N), T, d\right) = 
    \limsup_{\varepsilon\to 0}
   \frac{\sup_{f\in \mathcal{M}_1(\mathbb{C}P^N)}
   S\left(B_\delta\left(f, d_{\mathbb{C}}\right), \varepsilon\right)}{\log (1/\varepsilon)}. \]
So we need to study $B_\delta\left(f, d_{\mathbb{C}}\right)= \{g\in \mathcal{M}_1(\mathbb{C}P^N)|\, d_\mathbb{C}(f, g) \leq \delta\}$.

We sometimes need to consider balls in $\mathcal{M}_\lambda(\mathbb{C}P^N)$ for $\lambda\neq 1$.
So we introduce a new notation for clarifying the 
value of the parameter $\lambda$: 
For $\delta>0$ and $f\in \mathcal{M}_\lambda(\mathbb{C}P^N)$ we denote by 
$B_\delta\left(f, d_{\mathbb{C}}\right)_\lambda$ 
the $\delta$-ball around $f$ in $\mathcal{M}_\lambda(\mathbb{C}P^N)$ with respect to 
$d_{\mathbb{C}}$. Namely
\[ B_\delta\left(f, d_{\mathbb{C}}\right)_\lambda := \left\{g\in \mathcal{M}_\lambda(\mathbb{C}P^N)\middle|\, 
    \sup_{z\in \mathbb{C}} \mathbf{d}_{\mathrm{FS}}\left(f(z), g(z)\right) \leq \delta\right\}. \]
In this notation,
\[ \overline{\mdim}_M\left(\mathcal{M}_\lambda(\mathbb{C}P^N), T, d\right) = 
    \limsup_{\varepsilon\to 0} 
\frac{\sup_{f\in \mathcal{M}_\lambda(\mathbb{C}P^N)}S\left(B_\delta\left(f, d_{\mathbb{C}}\right)_\lambda, \varepsilon\right)}
{\log (1/\varepsilon)}. \]

\subsection{Nondegeneracy and key propositions} \label{subsection: nondegeneracy and key propositions}

The following definition is very important.
For $a\in \mathbb{C}$ and $R>0$ we denote $D_R(a) = \{z\in \mathbb{C}|\, |z-a|\leq R\}$.

\begin{definition}  \label{definition: nondegenerate Brody curves}
Let $f:\mathbb{C}\to \mathbb{C}P^N$ be a holomorphic curve.
For $R>0$ it is said to be \textbf{$R$-nondegenerate} if for all $a\in \mathbb{C}$
\begin{equation} \label{eq: nondegeneracy condition}
   \max_{z\in D_R(a)} |df|(z) \geq \frac{1}{R}. 
\end{equation}
\end{definition}

The basic idea behind this definition is as follows.
For $f\in \mathcal{M}_\lambda(\mathbb{C}P^N)$ we would like to study 
$B_\delta\left(f, d_{\mathbb{C}}\right)_\lambda$ by a deformation theory technique.
As is usual in deformation theory, the description of a small deformation of $f$ becomes simpler if 
a \textit{transversality condition} is satisfed at the curve $f$.
The above definition provides the transversality condition we need in our deformation theory\footnote{The condition 
in Definition \ref{definition: nondegenerate Brody curves} involves the parameter $R$.
So it is a kind of \textit{quantitative transversality condition}.}.
A trivial example which is \textit{not} $R$-nondegenerate is a constant curve (i.e. a holomorphic curve whose image is one-point).
Constant curves are the most ``singular'' object for our deformation theory.
The condition (\ref{eq: nondegeneracy condition}) means that $f$ is not close to a constant curve over the disk $D_R(a)$.
So an $R$-nondegenerate curve is not close to a constant curve over any $R$-disk\footnote{Yosida \cite{Yosida} defined that 
a meromorphic function $f\in \mathcal{M}_\lambda(\mathbb{C}P^1)$ is \textbf{of first category} if the closure of the $\mathbb{C}$-orbit
of $f$ does not contain a constant function. Definition \ref{definition: nondegenerate Brody curves} is a 
quantitative version of this old idea.}.
For example, nonconstant elliptic functions are $R$-nondegenerate for some $R>0$.

The next two propositions are slightly simpler versions of \cite[Proposition 3.2, Proposition 3.3]{Tsukamoto Brody}.
They are key results for the proof of (\ref{eq: upper bound on metric mean dimension}).

\begin{proposition} \label{prop: resolution}
There exist positive numbers $\delta_1$ and $C_1$ satisfying the following statement.
For any $\lambda>1$ there exists $R_1 = R_1(\lambda)>0$ such that 
for any $f\in \mathcal{M}_1(\mathbb{C}P^N)$ we can construct a map 
\[ \Phi: B_{\delta_1}\left(f, d_{\mathbb{C}}\right)_1 \to \mathcal{M}_\lambda(\mathbb{C}P^N) \]
satisfying 
\begin{enumerate}
   \item $\Phi(f)$ is $R_1$-nondegenerate.
   \item For any $g, h\in B_{\delta_1}\left(f,d_{\mathbb{C}}\right)_1$ and $z\in \mathbb{C}$ 
\begin{equation*}
  \begin{split}
   & \mathbf{d}_{\mathrm{FS}}\left(\Phi(g)(z), \Phi(h)(z)\right)  \leq C_1 \cdot \mathbf{d}_{\mathrm{FS}}\left(g(z), h(z)\right), \\
   & \mathbf{d}_{\mathrm{FS}}\left(g(z), h(z)\right)  \leq C_1 \cdot \max_{|w-z|\leq 3} \mathbf{d}_{\mathrm{FS}}\left(\Phi(g)(w), \Phi(h)(w)\right). 
  \end{split}
\end{equation*}
\end{enumerate}
\end{proposition}

This is a kind of ``resolution of singularity''. 
If we choose an arbitrary $f\in \mathcal{M}_1(\mathbb{C}P^N)$, it might be a degenerate curve.
However we can replace it by a nondegenerate one $\Phi(f)$.
The proof is an application of surgery (gluing).
Given $f\in \mathcal{M}_1(\mathbb{C}P^N)$, we look for its 
``degenerate region'' (i.e. the region where the norm $|df|$ is uniformly small).
We glue sufficiently many rational curves to $f$ over the degenerate region.
Then the resulting curve becomes nondegenerate.
The condition (2) means that this surgery procedure does not destroy the metric structure.

\begin{proposition}  \label{prop: deformation}
For any $R>0$ and $0<\varepsilon<1$ there exist positive numbers $\delta_2 = \delta_2(R)$, $C_2 = C_2(R)$ and 
$C_3 = C_3(\varepsilon)$ satisfying the following statement.
Let $f\in \mathcal{M}_2(\mathbb{C}P^N)$ be an $R$-nondegenerate curve, and let 
$\Lambda\subset \mathbb{C}$ be a square of side length $L\geq 1$. Then 
\[ \#\left(B_{\delta_2}\left(f, d_{\mathbb{C}}\right)_2, d_\Lambda, \varepsilon\right)
    \leq \left(\frac{C_2}{\varepsilon}\right)^{2(N+1)\int_{\Lambda} |df|^2\, dxdy + C_3\cdot L}. \]
\end{proposition}

This is proved by deformation theory. 
Given a nondegenerate curve $f$, we describe the ball $B_{\delta_2}\left(f, d_{\mathbb{C}}\right)_2$
 by a deformation theory technique and get the above estimate.

Here we do not provide the detailed proofs of the above two propositions.
They are the same as the proofs of \cite[Proposition 3.2, Proposition 3.3]{Tsukamoto Brody}.
(Indeed the above propositions are slightly simpler than \cite[Proposition 3.2, Proposition 3.3]{Tsukamoto Brody} 
because the argument of \cite{Tsukamoto Brody} is based on Lemma \ref{lemma: weaker version} whereas the 
argument here uses much simpler Theorem \ref{theorem: main result for R^D}.)
It is not our purpose here to explain the proofs of the above two propositions.
Our purpose is to explain how to use Theorem \ref{theorem: main result for R^D}.

\subsection{Proof of the upper bound} \label{subsection: proof of the upper bound}

Now we start to prove the upper bound (\ref{eq: upper bound on metric mean dimension}):
\begin{equation*} 
    \overline{\mdim}_M\left(\mathcal{M}_1(\mathbb{C}P^N), T, d\right) \leq 2(N+1) \rho_1(\mathbb{C}P^N). 
\end{equation*}
By Theorem \ref{theorem: main result for R^D}, for any $\delta>0$
\[ \overline{\mdim}_M\left(\mathcal{M}_1(\mathbb{C}P^N), T, d\right) = 
    \limsup_{\varepsilon\to 0} 
\frac{\sup_{f\in \mathcal{M}_1(\mathbb{C}P^N)}S\left(B_\delta\left(f, d_{\mathbb{C}}\right)_1, \varepsilon\right)}
{\log (1/\varepsilon)}. \]
Recall that 
\[ S\left(B_\delta\left(f, d_{\mathbb{C}}\right)_1, \varepsilon\right) = \limsup_{L\to \infty}
    \frac{\log\#\left(B_\delta\left(f, d_{\mathbb{C}}\right)_1, d_{[0,L]^2}, \varepsilon \right)}{L^2}. \]

Let $\delta_1$ and $C_1$ be positive constants introduced in Proposition \ref{prop: resolution}.
Take an arbitrary $1<\lambda<2$.
Let $R_1 = R_1(\lambda)$ be the positive constant introduced in Proposition \ref{prop: resolution} for this $\lambda$.
And, let $\delta_2 = \delta_2(R_1)$ be the positive constant introduced in 
Proposition \ref{prop: deformation} for this $R_1$.
Set 
\[ \delta = \min\left(\delta_1, \frac{\delta_2}{C_1}\right). \]

Take an arbitrary $f\in \mathcal{M}_1(\mathbb{C}P^N)$.
Applying Proposition \ref{prop: resolution} to $f$, there is a map 
$\Phi:B_\delta\left(f, d_{\mathbb{C}}\right)_1 \to \mathcal{M}_\lambda(\mathbb{C}P^N)$ such that 
\begin{itemize}
   \item $\Phi(f)$ is $R_1$-nondegenerate.
   \item For any $g, h\in B_{\delta}\left(f,d_{\mathbb{C}}\right)_1$ and $z\in \mathbb{C}$ 
\begin{equation}  \label{eq: Lipschitz}
    \mathbf{d}_{\mathrm{FS}}\left(\Phi(g)(z), \Phi(h)(z)\right)  
  \leq C_1 \cdot \mathbf{d}_{\mathrm{FS}}\left(g(z), h(z)\right), 
\end{equation}
\begin{equation} \label{eq: inverse Lipschitz}
   \mathbf{d}_{\mathrm{FS}}\left(g(z), h(z)\right) 
 \leq C_1 \cdot \max_{|w-z|\leq 3} \mathbf{d}_{\mathrm{FS}}\left(\Phi(g)(w), \Phi(h)(w)\right). 
\end{equation}
\end{itemize}
From (\ref{eq: Lipschitz}) and $\delta\leq \frac{\delta_2}{C_1}$,
for $g\in B_\delta(f, d_{\mathbb{C}})_1$
\[  \mathbf{d}_{\mathrm{FS}}\left(\Phi(f)(z), \Phi(g)(z)\right) \leq C_1\delta \leq \delta_2. \]
Hence
\[  \Phi\left(B_\delta\left(f, d_{\mathbb{C}}\right)_1\right) \subset 
      B_{\delta_2}\left(\Phi(f), d_{\mathbb{C}}\right)_\lambda  \subset B_{\delta_2}\left(\Phi(f),d_{\mathbb{C}}\right)_2. \]
From (\ref{eq: inverse Lipschitz}), for any $g, h\in B_\delta(f, d_{\mathbb{C}})_1$ and $L>0$
\[  d_{[0,L]^2}(g, h) \leq C_1 \cdot d_{[-3, L+3]^2}\left(\Phi(g), \Phi(h)\right). \]
Hence for any $0< \varepsilon<1$
\begin{equation*}
  \begin{split}
    \#\left(B_\delta\left(f, d_{\mathbb{C}}\right)_1, d_{[0,L]^2}, \varepsilon\right)
    & \leq \#\left(\Phi\left(B_\delta\left(f, d_{\mathbb{C}}\right)_1\right) , d_{[-3,L+3]^2}, 
              \frac{\varepsilon}{2C_1}\right) \\
    & \leq \#\left(B_{\delta_2}\left(\Phi(f), d_{\mathbb{C}}\right)_2 , d_{[-3,L+3]^2}, 
              \frac{\varepsilon}{2C_1}\right).
  \end{split}
\end{equation*}
We apply Proposition \ref{prop: deformation} to the curve $\Phi(f)$ and the square $[-3, L+3]^2$ of side length $L+6$.
Then
\[ \#\left(B_{\delta_2}\left(\Phi(f), d_{\mathbb{C}}\right)_2 , d_{[-3,L+3]^2}, \frac{\varepsilon}{2C_1}\right)  \leq 
    \left(\frac{2C_1 C_2}{\varepsilon}\right)^{2(N+1)\int_{[-3, L+3]^2} \left|d\Phi(f)\right|^2 \, dxdy + C_3\cdot (L+6)}. \]
Here $C_2 = C_2(R_1)$ and $C_3 = C_3\left(\frac{\varepsilon}{2C_1}\right)$ are the positive constants introduced in 
Proposition \ref{prop: deformation}. 
Therefore
\begin{equation*}
   \begin{split}
   & \log \#\left(B_\delta\left(f, d_{\mathbb{C}}\right)_1, d_{[0,L]^2}, \varepsilon\right)  \\ & \leq
   \log\left(\frac{2C_1 C_2}{\varepsilon}\right) \left\{2(N+1)\int_{[-3, L+3]^2} \left|d\Phi(f)\right|^2 \, dxdy + C_3\cdot (L+6)\right\}
   \end{split}
\end{equation*}
Divide this by $L^2$ and let $L\to \infty$. Then 
\[ S\left(B_\delta\left(f, d_{\mathbb{C}}\right)_1, \varepsilon\right) 
    \leq  \log\left(\frac{2C_1 C_2}{\varepsilon}\right) \left\{2(N+1)\rho\left(\Phi(f)\right)\right\}. \]
Noting $\Phi(f) \in \mathcal{M}_\lambda(\mathbb{C}P^N)$, we get 
\[ \sup_{f\in \mathcal{M}_1(\mathbb{C}P^N)}S\left(B_\delta\left(f, d_{\mathbb{C}}\right)_1, \varepsilon\right) 
     \leq 2(N+1) \rho_\lambda(\mathbb{C}P^N) \cdot \log\left(\frac{2C_1 C_2}{\varepsilon}\right). \]
Divide this by $\log(1/\varepsilon)$ and let $\varepsilon\to 0$.
Then (noting $C_1 C_2$ is independent of $\varepsilon$)
\begin{equation*}
   \begin{split}
   \overline{\mdim}_M\left(\mathcal{M}_1(\mathbb{C}P^N), T, d\right) & = 
    \limsup_{\varepsilon\to 0} 
\frac{\sup_{f\in \mathcal{M}_1(\mathbb{C}P^N)}S\left(B_\delta\left(f, d_{\mathbb{C}}\right)_1, \varepsilon\right)}{\log (1/\varepsilon)} \\
   &\leq 2(N+1) \rho_\lambda(\mathbb{C}P^N). 
   \end{split}
\end{equation*}
Recall that $1<\lambda<2$ is arbitrary and $\rho_{\lambda}(\mathbb{C} P^N) = \lambda^2 \rho(\mathbb{C}P^N)$.
So, letting $\lambda\to 1$, we get the result
\[   \overline{\mdim}_M\left(\mathcal{M}_1(\mathbb{C}P^N), T, d\right)  \leq 2(N+1)\rho_1(\mathbb{C}P^N). \]

The author recommends interested readers to compare the argument in this subsection with the one given in 
\cite[pp.947-949]{Tsukamoto Brody}.
Their basic structures are the same, but probably the argument here is a bit clearer.
This is because here we use Theorem \ref{theorem: main result for R^D} whereas \cite{Tsukamoto Brody} 
used Lemma \ref{lemma: weaker version}.

\vspace{0.5cm}

\address{ Masaki Tsukamoto \endgraf
Department of Mathematics, Kyushu University, Moto-oka 744, Nishi-ku, Fukuoka 819-0395, Japan}

\textit{E-mail}: \texttt{tsukamoto@math.kyushu-u.ac.jp}

\end{document}